# NETWORK ROUTING ON REGULAR DIRECTED GRAPHS FROM SPANNING FACTORIZATIONS

*Randall Dougherty*
*Center for Communications Research, Institute for Defense Analyses*
*La Jolla, California*

*Vance Faber*
*Center for Computing Sciences, Institute for Defense Analyses*
*Bowie, Maryland*

***Abstract.*** *Networks with a high degree of symmetry are useful models for parallel processor networks. In earlier papers, we defined several global communication tasks (universal exchange, universal broadcast, universal summation) that can be critical tasks when complex algorithms are mapped to parallel machines. We showed that utilizing the symmetry can make network optimization a tractable problem. In particular, we showed that Cayley graphs have the desirable property that certain routing schemes starting from a single node can be transferred to all nodes in a way that does not introduce conflicts. In this paper, we define the concept of spanning factorizations and show that this property can also be used to transfer routing schemes from a single node to all other nodes. We show that all Cayley graphs and many (perhaps all) vertex transitive graphs have spanning factorizations.*

**Introduction**. Networks with a high degree of symmetry are useful models for parallel processor networks. In earlier papers ([5], [8], [9]), we defined several global communication tasks (universal exchange, universal broadcast, universal summation) which can be critical tasks when complex algorithms are mapped to parallel machines. We showed that utilizing the symmetry can make network optimization a tractable problem. In particular, we showed in [9] (and earlier in [5]) that Cayley graphs have the desirable property that certain routing schemes starting from a single node can be transferred to all nodes in a way which does not introduce conflicts. In this paper, we extend this transference idea to a class of graphs that is more inclusive than Cayley graphs.

**Notation for graph theory**. This paper mainly focuses on directed graphs derived from groups. Here, a directed graph $G$ is a set of *vertices* $V$ and a collection $E$ of ordered pairs of distinct vertices $(u,v)$ called *edges*. We often let $n$ be the number of vertices and $m$ be the number of edges. If some pair appears more than once as an edge, then $G$ is called a *multigraph*. Otherwise the pairs form a set and $G$ is called *elementary*. The vertex $u$ in the pair is called the *tail* of the edge and the vertex $v$ is called the *head*.

*Definition* (Cayley coset graph). Let $\Gamma$ be a finite group, $H$ a subgroup and $\Delta$ a subset. Suppose

  (i)  $\Delta \cap H = \emptyset$ and $\Gamma$ us generated by $\Delta \cup H$,
  (ii) $H\Delta H \subseteq \Delta H$,
  (iii) $\Delta$ is a set of distinct coset representatives of $H$ in $\Gamma$.

Then we can form the *Cayley coset graph* $G = (\Gamma, \Delta, H)$ with the cosets $\{gH : g \in \Gamma\}$ as vertices and the set of pairs $(gH, g\delta H)$ with $\delta \in \Delta$ as edges. When $H$ is the identity subgroup, the graph is a *Cayley graph*.

A graph $G$ is *vertex transitive* if for any two vertices $u$ and $v$ there is an automorphism of $G$ which maps $u$ to $v$. The classic proof of Sabidussi [13] shows that a graph is vertex transitive if and only if it is a Cayley coset graph. An important aspect of the proof shows that one can construct a Cayley coset graph from a vertex transitive graph by using the automorphism group as the group $\Gamma$ required in the definition and the subgroup of automorphisms that fix a vertex as the required subgroup $H$. The generators $\Delta$ correspond to automorphisms that map a vertex to a neighbor..

**Spanning factorizations**. A 1-factor of a directed graph $G$ is a subgraph with both in-degree and out-degree equal to 1. (Some authors have called this a 2-factor. Our definition seems more consistent with the notation in undirected graphs. For example, if the edges are all bi-directional and the factor is a union of 2-cycles, then this would be an ordinary 1-factor in an undirected graph.) It is known that every regular directed graph with in-degree and out-degree $d$ has a 1-factoring with $d$ 1-factors. For completeness, we give the proof here.

*Fact 1. Every directed graph $G$ where the in-degree and out-degree of every vertex is $d$ has an edge disjoint decomposition into $d$ 1-factors.*

*Proof.* Form an auxiliary graph $B$ with two new vertices $u'$ and $u''$ for each vertex $u$. The edges of $B$ are the pairs $(u', v'')$ where $(u, v)$ is a directed edge in $G$. The undirected graph $B$ is bipartite and regular with degree $d$ and so by Hall's Marriage Theorem, it can be decomposed into $d$ 1-factors. Each of these 1-factors corresponds to a directed 1-factor in $G$.

In order to create a routing scheme for universal exchange (often called the transpose – see [10]) on $G$, we consider regular graphs with factorizations with additional properties.

*Definition.* Let $F_1$, $F_2$, $\cdots$, $F_d$ be the factors in a 1-factoring of $G$. We call a finite string of symbols from the set $\{F_i\}$ a *word*. If $v$ is a vertex and $\omega$ is a word, then $v\omega$ denotes the directed path (and its endpoint) in $G$ starting at $v$ and proceeding along the unique edge corresponding to each consecutive factor represented in the word $\omega$. If $G$ is a graph with $n$ vertices, we say that a 1-factoring and a set of $n$ words $W = \{\omega_0, \omega_1, \omega_2, \ldots, \omega_{n-1}\}$, $\omega_0 = \emptyset$, is a *spanning factorization* of $G$ (with word list $W$) if for every vertex $v$, the vertices $v\omega_i$ are distinct.

**Schedules**. A schedule for universal exchange associated with a factorization is an assignment of a time (a label) to each occurrence of each factor in the words of $W$ such that no time is assigned more than once to a particular factor and times assigned to the factors in a single word are increasing. The *time of a schedule* is the largest time assigned to any of the factors. If $T$ is the total time, the schedule can be thought of as a $d \times T$ array where each row corresponds to a factor and an entry in that row indicates which occurrence of that factor has been assigned the corresponding time. An entry in a row in the array can be empty indicating no occurrence of that factor has been assigned the given time. The power of a spanning factorization lies in the fact that a schedule can be used to describe an algorithm for conflict free global exchange of information between the vertices of the graph.

*Theorem 1. Suppose we have a schedule for a factorization of the graph $G$. Then the collection of directed label-increasing paths $v\omega_i$ for all $v$ and non-empty $\omega_i$ have the property that no edge in the graph is assigned the same time twice. A schedule for a spanning factorization yields a time labeled directed path between every two vertices so that no edge is labeled with the same time twice.*

*Proof.* Each edge in the graph is assigned to a single one factor. Assume there is an edge in the one factor $F$ that has been assigned the same time twice. Since every occurrence of $F$ in the words in $W$ has been assigned a unique time, this can only mean that there are two different vertices $u$ and $v$ and an initial subword $\omega$ of a word in $W$ such that the edges $(u\omega, u\omega F)$ and $(v\omega, v\omega F)$ are the same edge. Then $u\omega$ and $v\omega$ must be the same vertex. Let us assume that this is the shortest $\omega$ for which this happens. The word $\omega$ cannot be empty since $u$ and $v$ are different. But then the last factor in $\omega$ must also be the same edge, a contradiction. If we start with a spanning factorization, then all the non-empty paths from $v$ are unique, there are $n-1$ of them and none of them can return to $v$ so they must reach to every other vertex in the graph.

There are some additional properties that a spanning factorization with word list $W$ might have.

*Definition.* We say a spanning factorization is *balanced* if each factor appears nearly as often in the schedule as any other. We say the factorization is *short* if the average number of times a factor appears is the same as the theoretical lower bound $\theta$ based on the average distance between any two vertices and the number of edges. We say the factorization is *optimal* if it is short and balanced. A schedule $\Sigma$ is *minimum* for a spanning factorization, if it has time $\tau(\Sigma)$ equal to the theoretical minimum time for the factorization based on maximum number of times a factor appears. In mathematical terms, we can write

(i)      *balanced* if

$$\max_i |\{F_i : F_i \in W\}| = \left\lceil \sum_{i=1}^{d} |\{F_i : F_i \in W\}| / d \right\rceil;$$

(ii)      *short* if

$$\left\lceil \sum_{i=1}^{d} |\{F_i : F_i \in W\}| / d \right\rceil = \left\lceil \sum_{k=1}^{D} \frac{kN_k}{nd} \right\rceil,$$

(here $N_k$ is the number of times the distance between two vertices is $k$ and $D$ is the diameter);

(iii)      *optimal* if

$$\max_i |\{F_i : F_i \in W\}| = \left\lceil \sum_{k=1}^{D} \frac{kN_k}{nd} \right\rceil;$$

(iv)      *minimum* if

$$\tau(\Sigma) = \max_i |\{F_i : F_i \in W\}|.$$

Note that these parameters are ordered

$$\left\lceil \sum_{k=1}^{D} \frac{kN_k}{nd} \right\rceil \leq \left\lceil \sum_{i=1}^{d} |\{F_i : F_i \in W\}| / d \right\rceil \leq \max_i |\{F_i : F_i \in W\}| \leq \tau(\Sigma).$$

Creation of schedules for spanning factorizations are discussed in [11], where the following is proven.

*Fact 2. Every diameter two spanning factorization has a minimum schedule unless the max belongs to a factor $F_i$ which is not in a word of length one and is entirely absent in words of length two in one position, either first or second. In that case, the shortest time for any schedule is one more than the theoretical minimum.*

**Universal broadcast**. This paper concerns universal exchange (transpose). Employing these ideas for universal broadcast requires more restrictions on the list of words. In a universal broadcast, instead of sending a different piece of information to all other vertices, a vertex has one single piece of information to send to all others. To utilize this communication pattern, we impose an additional condition on the words in the list $W$. We say that the list $W$ is *hierarchical* if every initial subword of a word in $W$ is also a word in $W$. In addition, given a hierarchical list for universal broadcast, the list is thought of as a tree and each edge in the tree is labeled with a time only once. The problem of assigning an optimal schedule is greatly simplified because all that is needed is for the times on a factor to form a partial order. It still may be a difficult problem to find which is the best tree to use.

**Cayley coset graphs**. Our main goal is to find spanning factorizations for Cayley coset graphs. If the graph is a Cayley graph, this is easy.

*Theorem 2. Every Cayley graph has a short factorization.*

*Proof.* This is a sketch. Take a tree $T_1$ of shortest paths from the identity of the group. The factors consist of all the edges labeled with a specific generator. The words are just the paths in $T_1$, so the factorization is short.

*Question. Does every Cayley coset graph have a spanning factorization or even a short spanning factorization?*

**Example 1 – CP graphs**. CP graphs $G(d,D)$ are vertex symmetric digraphs with a large number of vertices for a given degree d and diameter D. They were first introduced in [6] and [7]. Many of the properties that make them desirable for multiprocessor networks have been studied in [1], [3] and [4]. In particular, [4] constructs broadcast trees for CP graphs which are related to our factorization. In [14], it is determined which CP graphs are Cayley graphs and thus these have a short factorization. We can show that all CP graphs have a short factorization.

*Theorem 3. $G(d,D)$ has a short factorization.*

*Proof.* There are several definitions of $G(d,D)$. First $G(d,D)$ is a Cayley coset graph of the symmetric group $\Gamma = S_{d+1}$ with generators $\Delta = \{(k...21) : 2 \leq k \leq d+1\}$ and coset subgroup $H$ the group fixing the elements $1, 2, \ldots, D$. $H$ is isomorphic to the symmetric group $S_{d+1-D}$. $G(d, D)$ can also be described as a graph with vertices that are vectors of length $D$ with distinct elements of $\{1, 2, \ldots, d+1\}$ as entries. To describe the edges, we use the notation in [2]. Let $x = x_1 x_2 \ldots x_D$ be an arbitrary vertex of $G$. We know that the outward neighbors of $x$ are of two types:

$$R_k(x_1 x_2 \ldots x_D) = x_k x_1 \ldots x_{k-1} x_{k+1} \ldots x_D \text{ for } 2 \leq k \leq D$$

$$S_m(x_1 x_2 \ldots x_D) = m x_1 x_2 \ldots x_{D-1} \text{ for } m \notin X = \{x_1, x_2, \ldots, x_D\}.$$

We let $F_1, F_2, \ldots, F_{D-1}$ be $R_2, R_3, \ldots, R_D$, that is, $F_j(x) = (x, R_{j+1}(x))$ for $j = 1, \ldots, D-1$. For $j \geq D$, we proceed as follows. First, consider the cyclic order on $\{1, 2, \ldots, d+1\}$. Index the members of $\overline{X}$, the complement of $X$, as $y_0, y_1, \ldots, y_{d-D}$ where $y_0$ is the first element of $\overline{X}$ following $x_D$. Then define $F_{j+D}(x) = (x, S_{y_j}(x))$ for $j = 0, 1, \ldots, d - D$. We claim $F_j$ defined in this way is a short factorization of $G(d, D)$. We have to prove two things. First that $\{F_j\}$ is a factorization and second that it is short. To this end, let $W$ be the set of words that express the unique shortest paths from $I = 12\ldots D$. These words are defined by the algorithm in Theorem 3.6 in [7]. Now consider the vertices given by $xw$, $w \in W$. Let $\alpha$ be the permutation defined by $\alpha(j) = x_j$ for $j = 1, 2, \ldots, D$ and $\alpha(j + D) = y_j$. Then if the word $w(F_1, \ldots, F_d)$ is a path from $I$ to $z_1 z_2 \ldots z_D$, then $xw$ is a path from $x$ to $\alpha(z_1)\alpha(z_2)\ldots\alpha(z_D)$. Since $\alpha$ induces an automorphism on $G$, all the $xw$, $w \in W$ are distinct. Furthermore, there are the right number of words of each length so these are the unique shortest paths starting from $x$. This shows we have a short factorization of $G(d, D)$.

*Theorem 4.* When $D = 2$ and $d > 2$ the factorization of Theorem 3 is both optimal and minimum. The theoretical minimum time cannot be achieved for $G(2, 2)$.

*Proof.* In [7], the number of vertices $n_k = N_k / n$ at distance $k$ from $I$ (or any other vertex) is given as

$$(d+1)_k - (d+1)_{k-1} = (d+1)_{k-1}(d+1-k)$$

where $(n)_k$ is the falling factorial. For $k = 2$, this is $d^2 - 1$. When the factor $F_1$ is repeated twice, it forms a directed 2-cycle so does not lead to a vertex with distance 2. All other pairs of factors are words in $W$. Thus each factor is used $2d$ times in a word of length 2, except $F_1$ which is only used $2d - 2$ times. Taking into account the words of length 1, we see that the time for the schedule is at least $2d + 1$. By Fact 2, the factorization of Theorem 3 has a minimum schedule with time $2d + 1$. The theoretical lower bound on the time is

$$\lceil (n_1 + 2n_2)/d \rceil = \lceil 1 + 2d - 2/d \rceil = 2d + 1$$

as long as $d > 2$. The graph $G(2,2)$ has 6 vertices, degree 2 and the sum of the distances from one vertex to all of the others is 8 so the theoretical minimum time is 4 not 5. The theoretical minimum cannot be achieved by using the unique shortest paths since they have 30 $F_2$ edges and at most 6 can be used at any one time. Replacing even one shortest path by a longer path will increase the lower bound on the time.

Now we can produce a minimum schedule for universal exchange using the factorization in Theorem 3 for $D > 3$. We start by looking at the usage of each $F_j$ in the factorization. To this end, we define a recursion that grows the tree of unique shortest paths from any vertex $v$.

*Theorem 5. Each vertex in the tree other than $v$ is assigned a pair of numbers $(c,t)$ with $c \geq 1$ and $t \geq 1$.*

*Initialization. The vertex $v$ is assigned $t = 0$. There are $d$ edges out from $v$, one for each of the factors $F_j$, $1 \leq j \leq d$. The vertices at the heads of these edges are assigned $c = j$ and $t = 1$.*

*Recursion. At a vertex with assignment $(c,t)$ with $t < D$, there is an out edge $F_j$ for any $j$ with $j \geq t$ except if $j = t$ and $c = 1$. If $j \geq t + c$ then the head of $F_j$ is assigned $(c, t+1)$. If $t \leq j < t + c$, then the head of $F_j$ is assigned $(c-1, t+1)$.*

*Proof.* This is just another expression of the algorithm in Theorem 3.5 in [7].

We often refer to first entry in $(c,t)$ as the $c$ label and the second as the $t$ label. The $t$ label denotes the distance from the vertex $v$ while the $c$ label can be thought of as keeping track of the number of remaining uses we are allowed for factors with small indices. Intuitively, these factors are a limited resource because they consist of short cycles and reusing them too often gives a path that does not increase distance to $v$.

*Lemma 6.* Suppose $t \leq D$. If the vertex $v$ is assigned $(c,t)$ then $c \leq d - t + 1$. In fact, there is a vertex assigned $(c,t)$ for every $c$ in the range $1 \leq c \leq d - t + 1$.

*Proof.* This is clearly true for $t = 1$ because $1 \leq c \leq d$. Also, moving out along the edges, the $c$ labels either remain the same or decrease by 1 and there is always at least one case where the $c$ label decreases. Suppose we consider the vertex labeled $(d - t + 1, t)$ and an edge belonging to the cycle $F_j$ going out to a vertex assigned $(d - t + 1, t + 1)$. Then we must have $j \geq t + (d - t + 1) = d + 1$ which is impossible.

*Question.* How many times does $F_j$ appear in paths of length $k$ for $k \leq D$?

Let $S_k(j,i)$ be the number of paths of length $k$ using $F_j$ as the $i$ th step. We can calculate $S(j,1) = S_k(j,1)$ as follows. First we calculate $T(c,t) = T_k(c,t)$, the number of leaves in the tree starting at $(c,t)$ and going out to distance $k$ from $v$. We know that $T(c,t)$ only has meaning if $1 \leq c \leq d - t + 1$; we define it to be zero elsewhere. Then $S(j,1) = T(j,1)$. Once we have $T$, we show a recursion in Lemma 13 that can be used to calculate $S$. We calculate $T(c,t)$ recursively.

*Lemma 7.* For $c = 1$,
$$T(1,t) = (d - t)T(1, t + 1) = (d - t)_{k-t}.$$

*Proof.* First, we calculate $T(1, k - 1)$. We are not allowed to reduce $c$ so we just have the $j$ with $j \geq t + 1 = k$. This gives $T(1, k - 1) = d - k + 1 = (d - t)_1$. Now we use induction. We compute $T(1, t - 1)$. Again we cannot reduce $c$ so we have the $j$ with $j \geq (t - 1) + 1 = t$. This gives

$$T(1, t - 1) = (d - t + 1)T(1,t) = (d - t + 1)(d - t)_{k-t}$$

by the induction hypothesis. But this is exactly $(d - t + 1)_{k-t+1}$.

*Lemma 8.* We have $T(c,k) = 1$ for $1 \leq c \leq d - k + 1$.

*Proof.* This is just a restatement of Lemma 6 and the fact that by definition, the tree starting at a vertex $v$ labeled $(c,k)$ has $v$ as its only leaf.

*Lemma 9.* For $1 \leq c \leq d - t + 1$,

$$T(c,t) = cT(c-1,t+1) + (d+1-t-c)T(c,t+1).$$

*Proof.* Note that $T(c,k) \leq 1$. When $c > 1$, there are two types of $j$: $t \leq j < t+c$ which gives $c$ vertices with label $(c-1, t+1)$ and leading to $cT(c-1, t+1)$ leaves; and $j \geq t+c$ which gives $d+1-t-c$ vertices with label $(c, t+1)$ and leading to $(d+1-t-c)T(c, t+1)$ leaves. When $c = 1$, the formula is given in Lemma 7.

Lemma 10. The value of $T(c,t)$ is given by:

(i) $T(c,t) = (d-t+1)_{k-t}$ when $k-t < c \leq d-t+1$;

(ii) $T(c,t) = (d-t+1)_{k-t} - (k-t)_c (d-t+1-c)_{k-t-c}$ when $1 \leq c \leq k-t$.

*Proof.* We prove this by induction, starting with $t = k$ and working backward.

For $t = k$, $k - t = 0$ so we are in case (i) and $(d-k+1)_0 = 1$ and thus this is a valid statement. Now we assume that $t < k$ and that both statements are valid for $t+1$. We break the proof that the statements are valid for $t$ up into cases.

First, suppose that $c = 1$. Then we have to verify statement (ii). We have

$$(d-t+1)_{k-t} - (k-t)_1(d-t)_{k-t-1} = (d-t+1)(d-t)_{k-t-1} - (k-t)(d-t)_{k-t-1} = (d-t)_{k-t}$$

which is known to be $T(1,t)$ by Lemma 7.

Second, suppose that $c > k - t \geq 1$. Then $c > c - 1 > k - 1 - t = k - (t+1)$ so both $T(c-1, t+1)$ and $T(c, t+1)$ are known from case (i) of the induction hypothesis. Using Lemma 9, we can compute

$$T(c,t) = cT(c-1,t+1) + (d+1-t-c)T(c,t+1)$$

$$= c(d-t)_{k-t-1} + (d+1-t-c)(d-t)_{k-t-1} = (d+1-t)(d-t)_{k-t-1} = (d+1-t)_{k-t}$$

as required.

Third, suppose that $c > 1$ and $c = k - t$. Then $c - 1 = k - (t+1)$ while $c > k - (t+1)$ and so $T(c-1, t+1)$ is known from case (ii) of the induction hypothesis and $T(c, t+1)$ is known from case (i). Using Lemma 9, we compute

$$T(c,t) = c[(d-t)_{k-t-1} - (k-t-1)_{c-1}(d-t+1-c)_{k-t-c}] + (d+1-t-c)(d-t)_{k-t-1}$$

$$= (d+1-t)(d-t)_{k-t-1} - c(k-t-1)_{c-1}(d-t+1-c)_{k-t-c}.$$

Since $c = k - t$, this yields $T(c,t) = (d+1-t)(d-t)_{k-t-1} - c!$ which matches statement (ii) in this case.

Finally, suppose that $c > 1$ and $c < k - t$ (so in addition, $t < k$). Then both $c$ and $c - 1$ are less than or equal to $k - (t+1)$ and $T(c-1, t+1)$ and $T(c, t+1)$ are known from case (ii) of the induction hypothesis. Using Lemma 9, we can compute

$$T(c,t) = cT(c-1, t+1) + (d+1-t-c)T(c, t+1)$$

$$= c[(d-t)_{k-t-1} - (k-t-1)_{c-1}(d-t-c+1)_{k-t-c}]$$

$$+ (d+1-t-c)[(d-t)_{k-t-1} - (k-t-1)_c(d-t-c)_{k-t-1-c}]$$

$$= (d+1-t)(d-t)_{k-t-1} - [c(k-t-1)_{c-1} + (k-t-1)_c](d-t-c+1)_{k-t-c}$$

$$= (d+1-t)(d-t)_{k-t-1} - (k-t)_c(d-t-c+1)_{k-t-c}$$

which matches statement (ii). This concludes the induction step and proves the Lemma.

Next we calculate two other auxiliary quantities. Let $V(c,t)$ be the number of vertices assigned the pair $(c,t)$.

*Lemma 11.* Let $V(c,t)$ be the number of vertices assigned the pair $(c,t)$ then

$$V(c, t+1) = (d+1)_t \text{ when } c \leq d - t$$

*and zero otherwise.*

*Proof.* By the initialization step in Theorem 5, $V(c,1) = 1$ for $1 \leq c \leq d$. Note that if $c \leq d - t$, then both $c \leq d - (t+1)$ and $c + 1 \leq d - (t+1)$. The recursion part of Theorem 5 then yields

$$V(c, t+1) = (d - t - c + 1)V(c,t) + (c+1)V(c+1, t).$$

Furthermore, this recursion and the initial condition are satisfied by $(d+1)_t$ which proves the Lemma.

*Lemma 12.* Let $U(j,c,t)$ be the total number of times that an edge in the factor $F_j$ appears in the tree starting at some vertex labeled $(c,t)$. Then

$U(j,c,t) = V(c,t)$ if $t < j \le d$ or $c \ne 1$ and $j = t$; otherwise $U(j,c,t) = 0$.

*Proof.* Since the presence of the factor $F_j$ depends only on the label $(c,t)$, $U(j,c,t)$ is either 0 or $V(c,t)$. The rules for growing the tree yield the result.

Now we deal with $S_k(j,t)$, $t \le k$, $1 \le j \le d$.

*Lemma 13.* The value of $S_k(j,t)$ can be calculated using $T(c,t)$:

$$S(j,1) = T(j,1)$$

and for all $j \ge t$

$$S_k(j,t+1) = (d+1)_{t-1}\left(\sum_{c=1}^{d-t} T(c,t+1) + T(j-t,t+1)\right)$$

where we define $T(0,t+1) = 0$.

*Proof.* First we mentioned above Lemma 7 that $S(j,1) = T(j,1)$. Now fix $j$, $c$ and $t$. Let $(c(j,c),t+1)$ be the label on the head of the edge in the factor $F_j$ starting at a vertex labeled $(c,t)$. Then $c(j,c) = c$ unless $t \le j < t+c$ in which case $c(j,c) = c-1$. The number of times that a single $F_j$ of this type is used by paths out to distance $k$ is then $T(c(j,c),t+1)$ and so edges of this type account for $U(j,c,t)T(c(j,c),t+1)$ appearances of $F_j$ at level $t+1$. Thus

$$S_k(j,t+1) = \sum_{c=1}^{d-t+1} U(j,c,t)T(c(j,c),t+1).$$

If $j > t$ then $U(j,c,t) = V(c,t)$ and

$$S_k(j,t+1) = \sum_{c=1}^{d-t+1} U(j,c,t)T(c(j,c),t+1) = \sum_{c=1}^{d-t+1} V(c,t)T(c(j,c),t+1)$$

$$= (d+1)_{t-1} \sum_{c=1}^{d-t+1} T(c(j,c), t+1) = (d+1)_{t-1} \left( \sum_{c=1}^{d-t} T(c, t+1) + T(j-t, t+1) \right)$$

If $j = t$ then $U(j,c,t) = V(c,t)$ for $c > 1$ and zero otherwise so

$$S_k(j, t+1) = \sum_{c=1}^{d-t+1} U(j,c,t) T(c(j,c), t+1) = \sum_{c=2}^{d-t+1} V(c,t) T(c-1, t+1)$$

$$= (d+1)_{t-1} \sum_{c=2}^{d-t+1} T(c-1, t+1) = (d+1)_{t-1} \sum_{c=1}^{d-t} T(c, t+1).$$

$$S_k(j,t) = (d+1)_{t-2}[(d-k+1)(d-t+2)_{k-t} + (d-t+1)_{k-t}$$
$$- (k-t)_{j-t+1}(d-j)_{k-j-1}],$$

*Theorem 14.* For all $t \leq k \leq D$, the number of paths of length $k$ using $F_j$ as the $t$ th step, $S_k(j,t)$, can be described (using the convention that the falling factorial is zero whenever the argument or the index is out of bounds) as follows:

1) if $t = 1$ then $S_k(j,1) = (d)_{k-1} - (k-1)_j (d-j)_{k-j-1}$,
2.1) if $t > 1$ and $j < t-1$ then $S_k(j,t) = 0$,
2.2) if $t > 1$ and $j \geq t-1$ then
$$S_k(j,t) = (d+1)_{t-2}[(d-k+1)(d-t+2)_{k-t} + (d-t+1)_{k-t}$$
$$- (k-t)_{j-t+1}(d-j)_{k-j-1}].$$

*Proof.* First, we expand the statement by eliminating the conventions. The equivalent statement has these cases

1) if $t = 1$ then

a) if $j \geq k$ then

$$S_k(j,1) = (d)_{k-1},$$

b) if $j < k$ then

$$S_k(j,1) = (d)_{k-1} - (k-1)_j (d-j)_{k-j-1};$$

2) if $t > 1$ then

a) if $j < t-1$ then

$$S_k(j,t) = 0,$$

b) if $j = t-1 = k-1$ then

$$S_k(j,t) = (d+1)_{t-2}(d-t+1),$$

c) if $j = t-1 < k-1$ then

$$S_k(j,t) = (d+1)_{t-2}(d-k+1)(d-t+2)_{k-t},$$

d) if $t-1 < j < k$ then

$$S_k(j,t) = (d+1)_{t-2}[(d-k+1)(d-t+2)_{k-t} + (d-t+1)_{k-t}$$

$$- (k-t)_{j-t+1}(d-j)_{k-j-1}],$$

e) if $j \geq k$ then

$$S_k(j,t) = (d+1)_{t-2}[(d-k+1)(d-t+2)_{k-t} + (d-t+1)_{k-t}].$$

Now we can prove each of the cases. We have seen that $S(j,1) = T(j,1)$. The two parts of the case $t=1$ are then given by evaluating the expression in Lemma 10 at $t=1$. For larger values of $t$, Theorem 5 only allows $F_j$ at the $t$th step if $j \geq t-1$. We can use Lemma 13 to calculate $S_k(j, t+1)$. First, we handle the summation. From Lemma 10 if $j = t = k-1$, then

$$\sum_{c=1}^{d-t} T(c, t+1) = \sum_{c=1}^{d-t} (d-t)_{k-t-1} = d-t.$$ This proves statement 2b.

Otherwise, if $k > t+1$ from Lemma 10 we have

$$\sum_{c=1}^{d-t} T(c, t+1) = \sum_{c=1}^{d-t} (d-t)_{k-t-1} - \sum_{c=1}^{k-t-1} (k-t-1)_c (d-t-c)_{k-t-1-c}$$

$$= (d-t)(d-t)_{k-t-1} - \sum_{c=1}^{k-t-1} (k-t-1)_c (d-t-c)_{k-t-1-c}.$$

Now we show that the summation on the far right satisfies

$$\sum_{c=1}^{k-t-1}(k-t-1)_c(d-t-c)_{k-t-1-c} = (k-t-1)(d-t)_{k-t-2}.$$

We use the well-known combinatorial identity

$$\sum_{b=1}^{p}\binom{a-b}{a-p} = \binom{a}{p-1}.$$

This has a falling factorial analog given by

$$\sum_{b=1}^{p}(p)_b(a-b)_{p-b} = p!\sum_{b=1}^{p}\binom{a-b}{a-p} = p!\binom{a}{p-1} = p(a)_{p-1}.$$

Setting $b=c$, $a=d-t$ and $p=k-t-1$ yields the desired identity:

$$\sum_{c=1}^{k-t-1}(k-t-1)_c(d-t-c)_{k-t-1-c} = (k-t-1)(d-t)_{k-t-2}.$$

This gives us

$$\sum_{c=1}^{d-t}T(c,t+1) = (d-t-1)(d-t)_{k-t-1} - (k-t-1)(d-t)_{k-t-2}$$

$$= \frac{(d-t)!}{(d-k+2)!}[(d-t-1)(d-k+2) - (k-t-1)]$$

$$= \frac{(d-t)!}{(d-k+2)!}(d-t+1)(d-k+1) = (d-k+1)(d-t+1)_{k-t-1}.$$

Finally, the value of $T(j-t,t+1)$ is obtained from Lemma 10. If $j=t$, $T(j-t,t+1)=T(0,t+1)=0$ by definition. This proves statement 2c. If $j>t$ and $k\leq j\leq d$ the value is $(d-t)_{k-t-1}$. This proves statement 2e. If $t<j<k$ the value has two terms

$$T(j-t,t+1) = (d-t)_{k-t-1} - (k-t-1)_{j-t}(d-j)_{k-j-1}..$$

This proves statement 2d.

*Lemma 15. At every level, cycles with larger indices appear more often than cycles with smaller indices, that is, for $t-1 \leq j < d$*

$$S_k(j+1,t) \geq S_k(j,t).$$

*Proof.* For $t=1$, this is clearly true by Theorem 14 part 1. If $t>1$, we notice that Lemma 13 gives

$$\Delta = S_k(j+1,t) - S_k(j,t) = T(j+1-t,t+1) - T(j-t,t+1).$$

As in Theorem 14, there are cases to consider from Lemma 10. When $j+1-t \geq k-(t+1)$ (that is, $j+1 \geq k-1$) the difference is clearly non-negative. When $1 \leq j < k-2$, the difference is

$$\Delta = -(k-t-1)_{j+1}(d-t-j-1)_{k-t-j} + (k-t-1)_j(d-t-j)_{k-t-j-1}$$

$$= \frac{(k-t-1)!(d-t-j-1)!}{(k-t-j-2)!(d-k-1)!}\left[-1 + \frac{d-t-j}{k-1-t-j}\right].$$

But $j < k-t-1$ and $k \leq d$ show the denominators are non-zero and $\Delta$ is positive.

*Theorem 16. There exists a minimum schedule $\Sigma$ for the factorization of Theorem 3.*

*Proof.* Let $\mu = \max_i |\{F_i : F_i \in W\}|$. By Lemma 15, $\mu = |\{F_d : F_d \in W\}|$, the number of times that $F_d$ appears on the words in $W$. We have seen by Theorem 5, that each vertex that is the tail of the edge corresponding to some occurrence of the factor $F_i$ is also the tail of an $F_d$. Assign $\mu$ distinct times to the $F_d$ in the words of $W$ so that the times are partially ordered by layer $t$. Now assign the same time to all factors occurring in $W$ with the same tail. This guarantees that no time will appear more than once on a given factor and that the times are partially ordered by layer.

*Discussion of the lack of balance in the transpose.* Let $\theta = \sum_{k=1}^{D} \frac{kn_k}{nd}$ and $\mu = |\{F_d : F_d \in W\}|$. By Theorem 14 part 2e,

$$\mu = \sum_{k=1}^{D}\sum_{t=1}^{k} S_k(d,t) = \sum_{k=1}^{D}\sum_{t=0}^{k-1}(d+1)_{t-2}[(d-k+1)(d-t+2)_{k-t} + (d-t+1)_{k-t}].$$

This expression can be manipulated to get

$$\mu = \sum_{k=1}^{D}(d+1)_k \sum_{t=0}^{k-1}\frac{1}{d-t+3}(\frac{d-k+1}{d-k+2} + \frac{1}{d-t+2}).$$

By the results in [7], we know that

$$\theta = \frac{1}{d}\sum_{k=1}^{D} k((d+1)_k - (d+1)_{k-1}) = \frac{1}{d}\sum_{k=1}^{D} k(d+1)_{k-1}(d-k+1).$$

It is clear from Lemma 15 that $\mu$ is diverging from $\theta$ as $D$ grows. Can we estimate this divergence?

**Example 2 – large graphs of diameter 2**. The vertex symmetric but non-Cayley graphs $H_q$ have been defined by McKay, Miller and Siran in [12]. They have $2q^2$ vertices, degree $(3q-1)/2$ and diameter 2, where $q = 4l+1$ is a prime power congruent to 1 (mod 4). They are the largest known vertex symmetric graphs of diameter 2 for a given degree. The vertices of $H_q$ have the form $(i,m,r)$ where $i,m \in F = GF(q)$ and $r \in Z_2$.

*Theorem 17. The graphs $H_q$ have a short spanning factorization.*

*Proof.* Let $z$ be a primitive root of $GF(q)$ and let $X = \{z^{2k}\}$ be the set of all even powers of $z$. Each vertex $(i,m,r)$ is adjacent to the $2l = (q-1)/2$ vertices $(i, m+z^r x, r)$ for $x \in X$ as well as to the $q$ vertices $(j, m+(-1)^r ij, 1-r)$ for $j \in F$. Thus the out-degree of $H_q$ is $d = (3q-1)/2$. We take as factors, the $(q-1)/2$ sets of edges $F_x = \{(i,m,r),(i,m+z^r x,r)\}$ for $x \in X$ and the $q$ sets of edges $F_j = \{(i,m,r),(i+j,m+(-1)^r i(i+j),1-r)\}$ $j \in F$. Note that $F_j$ followed by $F_{-j}$ takes any vertex back to itself. We call the first set of factors the *fix-r* factors and the second set the *cross-over* factors. Now we construct the set of words $W$. Each single factor is a word and $\emptyset$ is a word. There are three types of two-letter words. To describe the first, let $w$ be a square in $F$ such that $1+w$ is not a square. Then we take each word which is a product of the two fix-$r$ factors $F_x F_{xw}$; there are $(q-1)/2$ of these. For the second type, we use any two cross-over factors $F_{j_1}F_{j_2}$ except $F_jF_{-j}$; there are $q^2 - q$ of these. For the third type, we use one cross-over factor and one fix-$r$ factor in either order, $F_jF_x$ and $F_xF_j$; there are $2q(q-1)/2 = q^2 - q$ of these. Note that the number of these words is $2q^2$, the same as the number of vertices.

To show that this collection forms a spanning factorization, we need only show that all the $v\omega_i$ with $\omega_i \in W$ are different. To start out with, it is clear that words of the first and second type produce vertices which are distinct from the vertices produced by the third type. Also note that words of the first and second type also produce different vertices because the first type fixes the first coordinate while the second type does not. Some simple algebra shows that each of the members of each type produce vertices which are distinct from those produced by members of the same type. Finally, we have to show that no word with a single factor produces a vertex which is identical to that produced by a word with two factors. This is clearly true when comparing words of different parities, crossover and fix-$r$. Note that any single fix-$r$ factor preserves the first coordinate so it cannot match the result of two cross-over factors. Suppose there are two fix-$r$ factors $F_x$ and $F_y$ such that $vF_xF_{xw} = vF_y$. But this implies that

$$m + z^r x + z^r xw = m + z^r y,$$

so $x(1+w) = y$ which contradicts the fact that $1 + w$ is not a square. If $vF_xF_j = vF_k$, then $vF_x = vF_kF_{-j}$ and we have just ruled this out when $k$ and $j$ are unequal. If $k$ and $j$ are equal, this is not possible because $F_x$ does not fix $v$. If $vF_jF_x = vF_k$, then let $u = vF_j$ and then we have $uF_x = uF_{-j}F_k$ so this is also ruled out. This proves the first statement in the theorem. The second statement is true because the paths $v\omega_i$ are all of shortest possible length.

*Theorem 18.* The theoretical minimum time for transpose on $H_q$ is $\lceil 8q/3 \rceil$. However, the minimum time for the spanning factorization in Theorem 17 is $3q - 2$ and this is the best that can be done for any schedule.
*Proof.* For any graph of diameter two, the theoretical minimum time is

$$\left\lceil \frac{2(n-1)}{d} \right\rceil - 1 = \left\lceil \frac{8}{3}q - \frac{1}{9} + \frac{1}{9(3q-1)} \right\rceil.$$

In order to find the time for the minimum schedule, we have to calculate the number of times each factor appears. A fix-$r$ factor appears on the left in a word once with another fix-$r$ factor, $q$ times with a cross-over factor and once by itself.

A fix-$r$ factor appears on the right once with another fix-$r$ and $q$ times with a cross-over factor for a total of $2q+3$ times. A cross-over factor appears on the left in a word $q-1$ times with another cross-over factor, $(q-1)/2$ times with a fix-$r$ factor and once by itself. . A cross-over factor appears on the right in a word $q-1$ times with another cross-over factor and $(q-1)/2$ times with a fix-$r$ factor for a total of $3q-2$ times. Thus the minimum schedule has time $3q-2$. It can be achieved by Fact 2.

The only way to find a spanning factorization that has a smaller schedule would be to reduce the number of double cross-over factors. Potentially, we could do this by swapping double cross-over factors with some number of fix-$r$ factors. Unfortunately, that is not possible because every pair of cross-over factors changes the value of the first component in $(i,m,r)$ while fix-$r$ factors also fix the first component.

**The Cayley coset representation of $H_q$.** We can determine the Cayley coset representation of $H_q$, with the hope that another set of generators might produce a more balanced transpose. Let $B$ be the non-abelian group of order $q^3$ given by

$$B = <x, y, c : [x,y] = c, c \in Z(B), x^q = y^q = 1>$$

and $\eta$ be the generator of the cyclic group of order $2(q-1)$.

*Theorem 18. The graph $H_q$ can be represented as a Cayley coset graph of the semidirect product group $\Gamma = B \rtimes <\eta>$ with*

$$\eta c \eta^{-1} = c,$$
$$\eta x \eta^{-1} = xy^{-1}$$
$$\eta y \eta^{-1} = x^a y^{-1} c^{z/2}$$

*where $a = 1 + z$. The coset is given by $H = <\eta^2, \alpha>$ where $\alpha = x^a y^{-1} c c^{-az/2}$.*

*Proof.* First we examine some automorphisms of $H_q$. Let

$$f_s(i,m,r) = (i, m+s, r), s \in F$$
$$g_t(i,m,r) = (i+t, m - (-1)^r it + rt^2, r), t \in F$$

$$h(i,m,r) = ((-z)^r i, zm, 1-r).$$

Note that

$$h \circ g_t \circ h^{-1}(i,m,r) = (i + (-z)^{1-r} t, m - z^r it + z(1-r)t^2, r).$$

Also

$$h^2(i,m,r) = (-zi, z^2 m, r)$$
$$h^{-2}(i,m,r) = (-z^{-1} i, z^{-2} m, r).$$

With $\gamma = [g_1, h]$, we can calculate $[g_1, \gamma] = (f_1)^{-a}$. Some other useful relations are

$$f_s \circ f_t = f_{s+t}$$

$$g_s \circ g_t = f_{-ts} \circ g_{t+s}$$

$$h^2 g_t h^{-2} = g_{-zt}$$
$$(g_1)^k = f_{-k(k-1)/2} \circ g_k$$

$$h f_1 h^{-1} = f_1$$

$$h g_1 h^{-1} = g_1 \gamma.$$

Now we can identify $x = g_1$, $y = \gamma^{-1}$ and $c = (f_1)^{-a}$. We are going to show that we can identify $\eta = h$. To do this, we verify

$$hch^{-1} = h(f_1)^{-a} h^{-1} = (h f_1 h^{-1})^{-a} = (f_1)^{-a} = c$$

$$h x h^{-1} = h g_1 h^{-1} = g_1 \gamma = x y^{-1}.$$

It remains to check that $h y h^{-1} = x^a y^{-1} c^{z/2}$. We know from the relations above that

$$h^2 x h^{-2} = h^2 g_1 h^{-2} = g_{-z} = (f_1)^{z(z+1)/2} (g_1)^{-z} = (f_1)^{za/2} (g_1)^{-z} = c^{-z/2} x^{-z}.$$

Also

$$h^2 x h^{-2} = h(hxh^{-1})h^{-1} = h(xy^{-1})h^{-1} = hxh^{-1}hy^{-1}h^{-1} = xy^{-1}hy^{-1}h^{-1}.$$

Thus

$$c^{-z/2} x^{-z} = xy^{-1}hy^{-1}h^{-1}$$

so rearranging gives

$$h\, yh^{-1} = x^a y^{-1} c^{z/2}.$$

The group has order $2(q-1)q^3$ and it is easy to see that the subgroup $H$ has order $q(q-1)$ so the coset graph has the right number of vertices. We can also see that from the relations above that both $h^2$ and $\alpha = f_{1-a(a+1)/2}(g_1)^a \gamma$ fix the vertex $(0,0,0)$ of $H_q$. Since we have identified a group of automorphisms $\Gamma$ of $H_q$ and determined that the stabilizer of $(0,0,0)$ is the subgroup $H$, we have shown that the set of cosets of $H$ in $\Gamma$ along with the automorphisms that correspond to the edges of $H_q$ at $(0,0,0)$ must be the Cayley coset representation of $H_q$. We can give a set of distinct coset representatives for the edges at $(0,0,0)$. We claim that the $(q-1)/2$ elements $c^{-\beta/a}$ for $\beta = z^{2k}$ and the $q$ elements $y^j h^{-1}$ map $(0,0,0)$ to a distinct neighbor in $H_q$. This is easy to verify by direct computation using the definitions of $c$, $y$ and $h$.